\documentclass{amsart}
\usepackage{amsmath, amsfonts, amssymb}

\usepackage[all]{xy}
\SelectTips{cm}{10}\UseTips
\usepackage[all]{xy}
\usepackage{amsthm}
\SelectTips{cm}{10}\UseTips
\theoremstyle{plain}
\newtheorem{thm}{Theorem}
\newtheorem{prop}[thm]{Proposition}

\newtheorem{lemma}[thm]{Lemma}

\theoremstyle{definition}

\newtheorem{rem}[thm]{Remark}

\numberwithin{thm}{section}
\numberwithin{equation}{section}

\def\bbf{\mathbf}
\def\bfQ{\mathbf{Q}}
\def\a{\alpha}

\def\G{\Gamma}

\def\ov{\overline}
\def\c{\chi}

\def\quf{\mathbf{Q}}

\def\ideleF{\mathbf{A}_F^{\times}}

\def\adeleF{\mathbf{A}_F}
\def\adeleFf{\mathbf{A}_{F,f}}

\newcommand{\no}{\noindent}
\newcommand{\be}{\begin{equation}}
\newcommand{\ee}{\end{equation}}

\newcommand{\bmat}{\left[ \begin{matrix}}
\newcommand{\emat}{\end{matrix} \right]}
\newcommand{\bsmat}{\left[ \begin{smallmatrix}}
\newcommand{\esmat}{\end{smallmatrix} \right]}

\newcommand{\hs}{\hspace{2pt}}
\newcommand{\hf}{\hspace{5pt}}
\newcommand{\OF}{\mathcal{O}_F}

\newcommand{\OFq}{\mathcal{O}_{F,\mathfrak{q}}}

\newcommand{\Op}{\mathcal{O}_{(\mathfrak{p})}}
\newcommand{\fp}{\mathfrak{p}}
\newcommand{\fq}{\mathfrak{q}}
\newcommand{\fl}{\mathfrak{l}}

\newcommand{\fm}{\mathfrak{m}}
\newcommand{\fN}{\mathfrak{N}}
\newcommand{\fpi}{\mathfrak{p}_i}

\newcommand{\iy}{\infty}

\newcommand{\tm}{\tilde{M}}

\DeclareMathOperator{\Hom}{Hom}
\DeclareMathOperator{\Res}{Res}

\DeclareMathOperator{\Gal}{Gal}

\DeclareMathOperator{\Cl}{Cl}

\DeclareMathOperator{\restr}{res}
\DeclareMathOperator{\Frob}{Frob}
\DeclareMathOperator{\GL}{GL}
\DeclareMathOperator{\SL}{SL}

\begin{document}
\title{Ihara's lemma for imaginary quadratic fields}
\begin{abstract}

An analogue over imaginary quadratic fields of a result in algebraic 
number theory known as Ihara's lemma is established. More precisely, we 
show that for a prime ideal $\fp$ of the ring of integers of an 
imaginary quadratic field $F$, the kernel of the sum of the two standard 
$\fp$-degeneracy maps between the cuspidal sheaf cohomology $H^1_!(X_0, 
\tilde{M}_0)^2 \rightarrow H^1_!(X_1, \tilde{M}_1)$ is Eisenstein. Here 
$X_0$ and $X_1$ are analogues over $F$ of the modular curves $X_0(N)$ and 
$X_0(Np)$, respectively. To prove our theorem we use the method of modular 
symbols and the congruence subgroup property for the group $\SL_2$ which 
is due to Serre.

\end{abstract}

\author{Krzysztof Klosin}

\date{August 22, 2007}
\maketitle

\section{Introduction} \label{Introduction}

Ihara's lemma in the version stated in \cite{Ribet84} asserts that the 
kernel of the map $\a: J_0(N)^2 \rightarrow J_0(Np)$ is Eisenstein if 
$(N,p)=1$. Here $J_0(N')$ 
denotes the Jacobian of the compactified modular 
curve $\G_0(N')\setminus\mathbf{H}$, and $\a$ is the sum of the two 
standard $p$-degeneracy maps from $J_0(N)$ to $J_0(Np)$. The original 
proof of the result is due to Ihara \cite{Ihara75} and uses algebraic 
geometry. In \cite{Ribet84} Ribet gave a different proof without appealing 
to algebro-geometric methods.
The result was later improved upon by Khare \cite{Khare95} 
to dispose of the condition that $N$ be coprime to $p$. Khare also gives a 
rearranged proof in the case when $(N,p)=1$ using the method of modular 
symbols
(cf. \cite{Khare95}, Remark 
4). We will use his approach to generalize the result to imaginary 
quadratic fields, where 
algebro-geometric techniques are not available.

Let $F$ denote an imaginary 
quadratic extension of $\quf$ and 
$\OF$ its ring of integers. The reason why over $F$ the algebro-geometric 
machinery is 
not available is the fact that the symmetric space on which automorphic 
forms are 
defined is the hyperbolic-3-space, the product of $\bbf{C}$ and 
$\bbf{R}_+$, and the analogues $X_n$ of the modular curves are not 
algebraic varieties (cf. section \ref{Preliminaries}). However, 
\cite{Khare95} uses only group 
cohomology and his method may be adapted to the situation over an 
imaginary quadratic 
field. In this setting the Jacobians are replaced with certain sheaf 
cohomology groups $H^1_!(X_n, \tilde{M}_n)$ and for a prime $\fp \subset 
\OF$ we have analogues of the two standard $\fp$-degeneracy maps 
whose sum $H^1_!(X_0, \tilde{M}_0)^2 \rightarrow H^1_!(X_1, \tilde{M}_1)$ 
we will call $\a$. (For precise definitions see section 
\ref{Preliminaries}.)
The main result of this note (Theorem 
\ref{theoremmain}) then asserts that the 
kernel of $\a$ is Eisenstein (for definition of ``Eisenstein'' see section 
\ref{Main Result}). 

Originally Ihara's lemma had been used by Ribet \cite{Ribet84} to prove 
the existence of congruences between modular forms of level $N$ and those 
of level $Np$. His result, valid for forms of weight 2, was later 
generalized to arbitrary weight by Diamond \cite{Diamond91}, who used 
the language of cohomology like we chose to. A crucial ingredient in 
Diamond's proof is the self-duality of $H^1(\G_0(N), M)$. Over imaginary 
quadratic fields, as over $\bfQ$, there is a connection between the space 
of automorphic forms and the cohomology groups $H^1_!(X_n, \tilde{M}_n)$ 
called the Eichler-Shimura-Harder isomorphism (cf. \cite{Urban95}). 
However, there seems to be no obvious way to adapt the approach of Ribet 
and Diamond to our situation as $H^1_!(X_n, \tilde{M}_n)$ is not 
self-dual. 

Ihara's lemma was also used in the proof of modularity of Galois 
representations attached to elliptic curves over $\bfQ$ (\cite{Wiles95}, 
\cite{BCDT}). Thanks to the work of Taylor \cite{Taylor94} one can attach 
Galois representations to a certain class of automorphic forms on 
$\Res_{F/\bfQ} (\GL_{2 /F})$. One could hope that Ihara's lemma in our 
formulation could be useful in proving the converse to Taylor's theorem,
 i.e., that ordinary Galois 
representations of $\Gal(\ov{F}/F)$ (satisfying appropriate conditions)  
arise from automorphic forms, but at this moment this is a mere 
speculation as too many other important ingredients of a potential proof 
seem to be missing.

The author would like to thank Trevor Arnold, Tobias Berger, Brian Conrad, 
Chandrashekhar Khare and Chris 
Skinner for many helpful and inspiring discussions.

\section{Preliminaries}
\label{Preliminaries}

Let $F$ be an imaginary quadratic extension 
of $\quf$ and denote by $\OF$ its ring of integers. Let $\fN$ be an ideal 
of 
$\OF$ such that the $\mathbf{Z}$-ideal $\fN \cap \mathbf{Z}$ has a 
generator greater than 3. Let $\fp$ be a prime ideal such that $\fp 
\nmid\fN$. 
Denote 
by $\Cl_F$ the class group of $F$ and choose representatives of distinct 
ideal classes to be prime ideals $\fpi$, 
$i=1, \dots, \#\Cl_F$, relatively prime to both $\fN$ and $\fp$. Let 
$\pi$, 
(resp. $\pi_i$) be a uniformizer of the completion $F_{\fp}$ (resp. 
$F_{\fpi}$) of $F$ at the prime $\fp$ (resp. $\fpi$), and put 
$\tilde{\pi}$ (resp. $\tilde{\pi}_i$) to be the idele $(\dots, 1, \pi, 1, 
\dots) \in \ideleF$ (resp. $(\dots, 1, \pi_i, 1,
\dots) \in \ideleF$), where $\pi$ (resp. $\pi_i$) occurs at the $\fp$-th 
place 
(resp. $\fpi$-th place). We also put $\Op := \bigcup_{j=0}^{\iy} \fp^{-j} 
\OF$.

For each $n\in \mathbf{Z}_{\geq 0}$, we define compact open 
subgroups of $\GL_2(\adeleFf)$
$$K_n:= \left\{\bmat a&b \\ c&d \emat \in \prod_{\fq \nmid \iy} 
\GL_2(\OFq) \mid c \in \fN \fp^n \right\}.$$
Here $\adeleFf$ denotes the finite adeles of $F$ and $\OFq$ the 
ring of integers of $F_{\fq}$. For $n \geq 0$ we 
also set $K_n^{\fp} 
= \bmat \tilde{\pi}\\ &1 \emat K_n \bmat \tilde{\pi}\\ &1 \emat^{-1}$. 

For any compact open 
subgroup $K$ of $\GL_2(\adeleFf)$ we put $$X_K = \GL_2(F) \setminus 
\GL_2(\adeleF)/ K \cdot U_2(\mathbf{C}) \cdot Z_{\iy},$$ 
where $Z_{\iy} = \mathbf{C}^{\times}$ is the center of 
$\GL_2(\mathbf{C})$ and $U(2):= \{M \in \GL_2(\mathbf{C}) \mid M \ov{M}^t 
= 
I_2\}$ (here `bar' denotes complex conjugation and $I_2$ stands for the $2 
\times 
2$-identity matrix). If $K$ is sufficiently large (which will be the 
case for all compact open subgroups we will consider) this space is a 
disjoint union of $\#\Cl_F$ connected 
components $X_K = \coprod_{i=1}^{\# \Cl_F} (\Gamma_K)_i \setminus 
\mathcal{Z}$, where $\mathcal{Z} = \GL_2(\mathbf{C})/U_2(\mathbf{C}) 
\mathbf{C}^{\times}$ and $(\G_K)_i = \GL_2(F) \cap \bmat \tilde{\pi_i}\\ 
&1 
\emat 
K \bmat \tilde{\pi_i}\\ &1 \emat^{-1}$. To ease notation we put $X_n := 
X_{K_n}$, $X_n^{\fp} := X_{K_n^{\fp}}$, $\G_{n,i}:= (\G_{K_n})_i$ and 
$\G^{\fp}_{n,i}:= (\G_{K_n^{\fp}})_i$.

We have the following diagram:
\be \label{diagramK} \xy*!C\xybox{\xymatrix{\dots \ar[r]^{\subset} 
&K_{n+1}\ar[r]^{\subset}\ar[d]^{\wr}&
K_n \ar[r]^{\subset} \ar[d]^{\wr} & K_{n-1}\ar[r]^{\subset} \ar[d]^{\wr}& 
\dots\\
\dots \ar[r]^{\subset} &K^{\fp}_{n+1}\ar[r]^{\subset}\ar[ur]^{\subset}&
K^{\fp}_n \ar[r]^{\subset} \ar[ur]^{\subset} & 
K^{\fp}_{n-1}\ar[r]^{\subset} &
\dots
  }}\endxy \ee
where the horizontal and diagonal arrows are inclusions and the 
vertical arrows are conjugation by  the maps $\bmat 
\tilde{\pi}\\ 
&1\emat$. Diagram (\ref{diagramK}) is not commutative, but it is 
``vertically commutative'', by which we mean that given two objects in the 
diagram, two directed paths between those two objects define the same map 
if and only if the two paths contain the same number of vertical arrows.  

Diagram (\ref{diagramK}) induces the following vertically commutative 
diagram of the corresponding symmetric spaces:
\be \label{diagramX} \xy*!C\xybox{\xymatrix{\dots \ar[r]
&X_{n+1}\ar[r]\ar[d]^{\wr}&
X_n \ar[r] \ar[d]^{\wr} & X_{n-1}\ar[r] \ar[d]^{\wr}&
\dots\\
\dots \ar[r] &X^{\fp}_{n+1}\ar[r]\ar[ur]&
X^{\fp}_n \ar[r] \ar[ur] &
X^{\fp}_{n-1}\ar[r] &
\dots
  }}\endxy \ee
The horizontal and diagonal arrows in diagram (\ref{diagramX}) are the 
natural 
projections and the vertical arrows are maps given by $(g_{\iy}, g_f) 
\mapsto \left(g_{\iy}, g_f \bsmat \tilde{\pi}\\
&1\esmat^{-1} \right)$. 

Let $M$ be a torsion abelian group of exponent relatively prime to 
$\#\OF^{\times}$ endowed with a $\GL_2(F)$-action. Denote by $\tilde{M}_K$ 
the sheaf of continuous sections of the 
topological covering $\GL_2(F)\setminus (\GL_2(\adeleF)/K 
\cdot U_2(\mathbf{C}) \cdot 
Z_{\iy})\times M \rightarrow X_K$, where $\GL_2(F)$ acts diagonally on 
$(\GL_2(\adeleF)/K U_2(\mathbf{C})
Z_{\iy})\times M$. Here $M$ is equipped with the discrete topology. Since 
we 
will only be concerned with the case when $M$ is a trivial 
$\GL_2(F)$-module, we assume it from now on. This means that $\tilde{M}_K$ 
is a constant sheaf. As above, we put $\tilde{M}_n:= \tilde{M}_{K_n}$ and 
$\tilde{M}^{\fp}_n:=\tilde{M}_{K^{\fp}_n}$.

Given a surjective map $\phi: X_K \rightarrow X_{K'}$, we get an 
isomorphism of 
sheaves 
$\phi^{-1}\tilde{M}_{K'} \xrightarrow{\sim} \tilde{M}_K$, which yields a 
map on 
cohomology $$H^q(X_{K'}, \tilde{M}_{K'}) \rightarrow H^q(X_{K}, 
\phi^{-1}\tilde{M}_{K'}) \cong H^q(X_K, \tilde{M}_{K}).$$ Hence diagram 
(\ref{diagramX}) gives rise to a vertically commutative diagram of 
cohomology groups:
\be \label{diagramc} \xymatrix@C4em@R5em{\dots 
&H^q(X_{n+1}, \tilde{M}_{n+1})\ar[l]&
H^q(X_n,\tilde{M}_n) \ar[l]_{\a_1^{n,n+1}} \ar[dl]_{\a_1^{n, (n+1)\fp}}& 
H^q(X_{n-1}, \tilde{M}_{n-1})\ar[l]_{\a_1^{n-1,n}}\ar[dl]_{\a_1^{n-1, 
n\fp}}&
\dots\ar[l]\\
\dots &H^q(X^{\fp}_{n+1}, 
\tilde{M}^{\fp}_{n+1})\ar[l]\ar[u]_{\wr}^{\a_{\fp}^{n+1}}&
H^q(X^{\fp}_n, \tilde{M}^{\fp}_n) \ar[l]_{\a_1^{n\fp, (n+1)\fp}} 
\ar[u]_{\wr}^{\a_{\fp}^{n}} &
H^q(X^{\fp}_{n-1}, \tilde{M}^{\fp}_{n-1})\ar[l]_{\a_1^{(n-1)\fp, n\fp}} 
\ar[u]_{\wr}^{\a_{\fp}^{n-1}}&
\dots\ar[l]
  } \ee
These sheaf cohomology groups can be related to the group cohomology 
of $\G_{n,i}$ and $\G_{n,i}^{\fp}$ with coefficients in $M$. In fact, for 
each compact open 
subgroup $K$ with corresponding decomposition $X_K = \coprod_{i=1}^{\# 
\Cl_F} 
(\Gamma_K)_i \setminus
\mathcal{Z}$, we have the following 
commutative diagram in which the horizontal maps are inclusions:
\be\label{diagramgc}\xy*!C\xybox{\xymatrix{H^q_!(X_K, \tilde{M}_K) \ar[r]  
&H^q(X_K, 
\tilde{M}_K) \\
\bigoplus_{i=1}^{\#\Cl_F} H^q_P((\G_K)_i, M) \ar[r] \ar[u]& 
\bigoplus_{i=1}^{\#\Cl_F} H^q((\G_K)_i, M) \ar[u]}}\endxy\ee 
Here $H^q_!(X_K, \tilde{M}_K)$ denotes the image of the cohomology 
with compact support $H^q_c(X_K, \tilde{M}_K)$ inside $H^q(X_K, 
\tilde{M}_K)$ and $H^q_P$ denotes the parabolic cohomology, i.e., 
$H^q_P((\G_K)_i, M):= \ker(H^q((\G_K)_i, M)\rightarrow \bigoplus_{B \in 
\mathcal{B}} H^q((\G_K)_{i,B}, M))$, where $\mathcal{B}$ is the set of 
Borel subgroups of 
$GL_2(F)$ and $(\G_K)_{i,B}:= (\G_K)_{i} \cap B$. The vertical 
arrows in diagram (\ref{diagramgc}) are isomorphisms provided that there 
exists a torsion-free normal 
subgroup of $(\G_K)_i$ of finite index relatively prime to 
the exponent of $M$. If $K=K_n$ or $K=K_n^{\fp}$, $n\geq 0$, this 
condition is satisfied because of our assumption that $\fN \cap 
\mathbf{Z}$ has a generator 
greater than 3 and the exponent of $M$ is relatively prime to $\# 
\OF^{\times}$ (cf. \cite{Urban95}, 
section 2.3). In what follows 
we may therefore identify the sheaf cohomology with the group cohomology. 
Note that all maps in 
diagram (\ref{diagramc}) preserve parabolic cohomology. The maps 
$\a_1^{*,*}$ are the natural restriction maps on group cohomology, so 
in 
particular they preserve the decomposition $\bigoplus_{i=1}^{\#\Cl_F} 
H^q((\G_K)_i, M)$.  

Using the identifications of diagram (\ref{diagramgc}) we can prove the 
following result which will be useful later:

\begin{lemma} \label{lemmainjective}

The map $\a_1^{0,1}: H^1_!(X_0, \tilde{M}_0) \rightarrow 
H^1_!(X_{1}, \tilde{M}_{1})$ is injective.

\end{lemma}

\begin{proof} Using the isomorphism between group and sheaf cohomology 
all we 
need to prove is that the restriction maps $\restr_i: H^1(\G_{0,i},M) 
\rightarrow H^1(\G_{1,i},M)$ are injective. Since $M$ is a trivial 
$\G_{0,i}$-module, the cohomology groups are just $\Hom$s, so it is 
enough to show the following statement: if $G$ denotes the smallest normal 
subgroup of 
$\G_{0,i}$ 
containing 
$\G_{1,i}$, then $G=\G_{0,i}$. For this we use the decomposition 
$$\G_{0,i} 
=\coprod_{\substack{k\in R(\OF/\fp) \\ k \in \fN }} \G_{1,i} \bmat 
1 \\ k&1 \emat \sqcup \G_{1,i} \bmat A&B \\ C&D \emat,$$
where the 
matrix $\bsmat A&B\\C&D\esmat$ is chosen so that $C$ 
and $D$ are relatively prime elements of $\OF$ with $C\in \fN$,
 $D\in \fp$, and $A \in \OF$, $B \in \fp_i$ satisfy $AD-BC=1$. Here 
$R(\OF/\fp)$ denotes a set of representatives in $\OF$ of 
the distinct residue classes of $\OF/\fp$. 
Let $\bsmat a&b\\c&d \esmat\in \G_{1,i}$ with $d\in \fp$. Then for any $k 
\in 
\fN \fp_i^{-1}\OF$ we have 
$$\bmat a&b \\c&d \emat^{-1} \bmat 1+bdk& -b^2k 
\\ d^2  k& 
1-bdk \emat \bmat a&b \\c&d \emat = \bmat 1\\k & 1 \emat$$
and the matrix 
$\bsmat 1+bdk& -b^2k \\ d^2  k&
1-bdk \esmat \in \G_{i,1}$, hence $G$ contains $\coprod_{\substack{k\in 
R(\OF/\fp) \\ k \in \fN }} \G_{1,i} \bsmat
1 \\ k&1 \esmat$, and thus $G=\G_{i,0}$.  \end{proof}

We can augment diagram (\ref{diagramK}) on the right by introducing one 
more 
group:
$$K_{-1}:= \left\{\bmat a&b \\ c&d \emat \in \GL_2(F_{\fp})\times 
\prod_{\fq
\nmid \fp \iy}
\GL_2(\OFq) \mid c \in \fN , \hs ad-bc \in \prod_{\fq \nmid \iy}
\OFq^{\times}
\right\}.$$
The group $K_{-1}$ is not compact, but we can still define 
$$\G_{-1,i}:= 
\GL_2(F) \cap \bmat \tilde{\pi}_i \\ &1 \emat K_{-1} \bmat \tilde{\pi}_i 
\\ 
&1 \emat^{-1}$$
for $i=1, \dots \#\Cl_F$. After identifying 
the sheaf cohomology 
groups $H^1_!(X_0, \tilde{M}_0)$ and $H^1_!(X_0^{\fp}, \tilde{M}_0^{\fp})$ 
with the groups $\bigoplus_i H^1_P(\G_{0,i}, M)$ and $\bigoplus_i 
H^1_P(\G^{\fp}_{0,i}, M)$, respectively, using diagram (\ref{diagramgc}), 
we 
can augment diagram (\ref{diagramc}) on the right in the following way
\be \label{diagramcaug} \xy*!C\xybox{\xymatrix@C5em@R7em{\dots
&H^q_!(X_{1}, \tilde{M}_{1})\ar[l]&
H^q_!(X_0,\tilde{M}_0) \ar[l]_{\a_1^{0,1}} \ar[dl]_{\a_1^{0, 1\fp}}&
H^q_!(X_{-1}, \tilde{M}_{-1})\ar[l]_{\a_1^{-1,0}}\ar[dl]_{\a_1^{-1,
0\fp}}\\
\dots &H^q_!(X^{\fp}_{1},
\tilde{M}^{\fp}_{1})\ar[l]\ar[u]_{\wr}^{\a_{\fp}^{1}}&
H^q_!(X^{\fp}_0, \tilde{M}^{\fp}_0) \ar[l]_{\a_1^{0\fp, 1\fp}}
\ar[u]_{\wr}^{\a_{\fp}^{0}}}}\endxy \ee
Here we put $H^q_!(X_{-1}, \tilde{M}_{-1}):=\bigoplus_i H^1_P (\G_{-1,i}, 
M)$ and the maps $\a_{1}^{-1,0}$ and $\a_{_1}^{-1,0\fp}$ are direct sums of the 
restriction maps.

The sheaf and group cohomologies are in a natural way modules over the 
corresponding Hecke 
algebras. (For the definition of the Hecke action on cohomology, see 
\cite{Urban95} or \cite{Hida93}). Here we will only consider the 
subalgebra 
$\mathbf{T}$ of the full Hecke algebra which is generated over $\bbf{Z}$ 
by the 
double cosets $T_{\fp'}:=K \bmat\pi'\\
&1\emat K$ and $T_{\fp',\fp'}:=K \bmat \pi' \\ &\pi' \emat K$ for $\pi'$ 
a uniformizer of $F_{\fp'}$ with $\fp'$ running over prime ideals of $\OF$ 
such that 
$\fp' \nmid \fN \fp$. The algebra $\bbf{T}$ acts on all the cohomology 
groups in diagram (\ref{diagramcaug}). 
Moreover if $\fp'$ is principal, the induced action of $T_{\fp'}$ and 
$T_{\fp', \fp'}$ 
on 
the group cohomology respects the decomposition $\bigoplus_i 
H^q_*((\G_K)_i, 
M)$, where 
$*=\emptyset$ or $P$.

\section{Main result} \label{Main Result}

We will say that a maximal ideal $\mathfrak{n}$ of the Hecke algebra 
$\mathbf{T}$ is \textit{Eisenstein} if $T_{\mathfrak{l}}\equiv 
N\mathfrak{l} +1\hf 
(\textup{mod} \hs \mathfrak{n})$ for all ideals $\mathfrak{l}$ of $\OF$ 
which 
are trivial as elements of the ray class group of conductor 
$\mathfrak{n}$. Such 
ideals $\fl$ are principal and have a generator $l$ with 
$l \equiv 1 \hf (\textup{mod} \hs \mathfrak{n})$. Here $N\mathfrak{l}$ 
denotes 
the ideal norm. 

From now on we fix a non-Eisenstein maximal ideal $\mathfrak{m}$ of the 
Hecke 
algebra $\mathbf{T}$. Our main result is the following theorem.

\begin{thm} \label{theoremmain} Consider the map 
$H^1_!(X_0,\tilde{M}_0)^2 
\xrightarrow{\alpha} H^1_!(X_1, \tm_1)$ defined as $\alpha: (f,g) \mapsto 
\a_1^{0,1}f+\a_{\fp}^1 \a_1^{0, 1\fp}g$. The localization $\a_{\fm}$ of 
$\a$ is injective. \end{thm}

We prove Theorem \ref{theoremmain} in two steps. Define a map 
$\beta:H^1_!(X_{-1}, \tm_{-1}) \rightarrow H^1_!(X_0,\tilde{M}_0)^2$ by 
$g' \mapsto (-\a_{\fp}^0\a_1^{-1,0\fp}g', \a_1^{-1,0} g')$ and note that 
$\alpha \beta =0$ by the vertical commutativity of diagram 
(\ref{diagramcaug}), i.e., $\ker \a \supset \beta \bigl(
H^1_!(X_{-1}, \tm_{-1})\bigr)$. We first 
prove

\begin{prop} \label{theorem2} $\ker \a = \beta \bigl( 
H^1_!(X_{-1}, \tm_{-1})\bigr)$. \end{prop}

\no Then we show 

\begin{prop} \label{theorem3} $H^1_!(X_{-1}, \tm_{-1})_{\fm}=0$. 
\end{prop}

\no Propositions \ref{theorem2} and \ref{theorem3} imply 
Theorem \ref{theoremmain}.

The idea of the proof is due to Khare \cite{Khare95} and uses 
modular 
symbols, which we now define. Let $D$ denote the free abelian 
group on the set $\mathcal{B}$ of all Borel subgroups of $\GL_2(F)$. The 
action of 
$\GL_2(F)$ on $\mathcal{B}$ by conjugation gives rise to a 
$\mathbf{Z}$-linear action of $\GL_2(F)$ on $D$. We sometimes identify 
$\mathcal{B}$ 
with $\mathbf{P}^1(F)=\bigl\{\frac{a}{c}\mid a\in \OF,c\in \OF\setminus 
\{0\}\bigr\}\cup \{\iy\}$, 
on which $\GL_2(F)$ acts by the linear fractional 
transformations. 
Let $D_0:=\{\sum n_i B_i \in D \mid n_i \in \mathbf{Z}, B_i \in 
\mathcal{B}, \sum n_i=0\}$ be the subset of elements of degree zero. If 
$X_K = \coprod_i \G_i\setminus\mathcal{Z}$, then 
for each $\G_i$ the exact sequence
$$0 \rightarrow \Hom_{\mathbf{Z}}(\mathbf{Z}, M) 
\rightarrow 
\Hom_{\mathbf{Z}}(D,M) \rightarrow \Hom_{\mathbf{Z}}(D_0,M) \rightarrow 
0$$
gives rise to an exact sequence
\begin{multline} \label{exactsequence}0 \rightarrow M \rightarrow 
\Hom_{\mathbf{Z}[\G_i]}(D,M) \rightarrow 
\Hom_{\mathbf{Z}[\G_i]} (D_0,M) \rightarrow \\
\rightarrow H^1(\G_i,M) \rightarrow  
H^1(\G_i, \Hom_{\mathbf{Z}}(D,M)).\end{multline}
The group $\Hom_{\mathbf{Z}[\G_i]} (D_0,M)$ is called the group of 
\textit{modular symbols}. 

\begin{lemma} \label{shapirolemma} Let $\G$ be a group acting on the set 
$\mathcal{B}$ of Borel subgroups of $\GL_2(F)$ and let $C$ denote a set 
representatives for the $\G$-orbits of $\mathcal{B}$. Then for any trivial 
$\G$-module $W$,
$$H^1(\G, \Hom_{\mathbf{Z}}(D,W)) = \bigoplus_{c\in C} H^1(\G_c,W),$$
where $\G_c$ is the stabilizer of $c$ in 
$\G$. \end{lemma}

\begin{proof} The $\Gamma$-module structure on
$\text{Hom}_{\mathbf{Z}}(D,W)$ is defined via $\phi ^{\gamma} (x) = \phi
(\gamma^{-1} x)$ and on $ \bigoplus_{c\in C}
\text{Ind}_{\Gamma_c}^{\Gamma} W$ via $$(f_{c_1},
\dots, f_{c_n})^{\gamma} (\gamma_{c_1},\dots, \gamma_{c_n}) = (f_{c_1},
\dots, f_{c_n})(\gamma_{c_1}\gamma,\dots, \gamma_{c_n}\gamma).$$
Note that we
have a $\Gamma$-module isomorphism $\Phi 
:\text{Hom}_{\mathbf{Z}}(D,W)\xrightarrow{\sim} \bigoplus_{c\in C}
\text{Ind}_{\Gamma_c}^{\Gamma} W$, given by $\Phi (\phi)_c (\gamma)=\phi
(\gamma^{-1}c)$.
Thus $$H^1(\Gamma,\text{Hom}_{\mathbf{Z}}(D,W)) \simeq H^1\bigl(\Gamma,
\bigoplus_{c\in C} \text{Ind}_{\Gamma_c}^{\Gamma} W\bigr) \simeq 
\bigoplus_{c
\in C} H^1(\Gamma,\text{Ind}_{\Gamma_c}^{\Gamma} W)$$ since the action of
$\Gamma$ stabilizes $\text{Ind}_{\Gamma_c}^{\Gamma} W$ for every $c \in 
C$.
The last group is in turn isomorphic to $\bigoplus_{c \in C}
H^1(\Gamma_c,W)$ by Shapiro's Lemma. \end{proof}

By taking the direct sum of the exact sequences 
(\ref{exactsequence}) and using Lemma \ref{shapirolemma}, we obtain the 
exact sequence
\begin{multline}  \label{exactsequence15} 0 \rightarrow 
\bigoplus_i M 
\rightarrow
\bigoplus_i\Hom_{\mathbf{Z}[\G_i]}(D,M) \rightarrow\\
\rightarrow \bigoplus_i\Hom_{\mathbf{Z}[\G_i]} (D_0,M) 
\rightarrow \bigoplus_i 
H^1_P(\G_i,M) 
\rightarrow 
0,\end{multline}
where the last group is isomorphic to $H^1_!(X_K, \tilde{M}_K)$. 

\begin{rem} The space of modular symbols 
$\bigoplus_i\Hom_{\mathbf{Z}[\G_i]} (D_0,M)$ is also a Hecke module in a 
natural way. In fact it can be shown (at least if $\fN$ is 
square-free) that the localized map 
$\bigl(\bigoplus_i\Hom_{\mathbf{Z}[\G_i]}(D,M)\bigr)_{\fm} \rightarrow
\bigl(\bigoplus_i\Hom_{\mathbf{Z}[\G_i]} (D_0,M)\bigr)_{\fm}$ is an 
isomorphism, but we will not need this fact.

\end{rem}

\section {Proof of Proposition \ref{theorem2}} \label{Proof1}

Suppose $(f,g) 
\in \ker \a$, i.e., $\a_1^{0,1} f = -\a_{\fp}^1 \a_1^{0,1\fp} g$. Let $h_1 
= -g \in H^1(X_0, \tilde{M}_0)$ and $h_2\in H^1(X_0^{\fp}, 
\tilde{M}_0^{\fp})$ be the pre-image of $f$ under the isomorphism 
$\a_0^{\fp}$. Then $\a_1^{0,1} 
\a_{\fp}^0 h_2 =\a_{\fp}^1 \a_1^{0,1\fp} h_1$. By the vertical 
commutativity of diagram (\ref{diagramc}), we have $\a_1^{0,1} \a_{\fp}^0 
= \a_{\fp}^1 \a_1^{0\fp, 1\fp}$, whence $\a_{\fp}^1 \a_1^{0\fp, 1\fp} h_2 
= 
\a_{\fp}^1 \a_1^{0,1\fp} h_1$. Since $\a_{\fp}^1$ is an isomorphism, we 
get $\a_1^{0\fp, 1\fp} h_2 = \a_1^{0,1\fp} h_1 \in H^1(X_1^{\fp}, 
\tilde{M}^{\fp}_1)$. 

For $K\subset K'$ two compact open subgroups of $\GL_2(\adeleFf)$ for 
which $X_K = \coprod_i \G_i\setminus \mathcal{Z}$ and $X_{K'} = 
\coprod_i \G'_i\setminus \mathcal{Z}$ with $\G_i, 
\G'_i \in \mathcal{G}_i$, $\G_i \subset \G'_i$, we have a commutative 
diagram
\be
\label{diagraminclusion0}\xymatrix{\Hom_{\mathbf{Z}[\G_i]}(D_0,M)
\ar[r]^-{\phi_{\G_i}}& H_P^1(\G_i, M)\\
\Hom_{\mathbf{Z}[\G'_i]}(D_0,M)
\ar[r]^-{\phi_{\G'_i}}\ar[u]^{\textup{inclusion}}& 
H_P^1(\G'_i,
M)\ar[u]^{\textup{res}}}, \ee
where the maps $\phi_{\G_i}$ and $\phi_{\G'_i}$ denote the appropriate 
connecting homomorphisms from exact sequence (\ref{exactsequence15}). So 
far
we have shown that
\be\label{sofar} -\a^{0,1\fp}_1 g = \a_1^{0,1\fp}h_1 = \a_1^{0\fp, 1\fp}
h_2.\ee
We identify $g$ with  a tuple $(g_i)_i \in \bigoplus_i H^1_P(\G_{0,i},M)$ 
and define $(h_1)_i\in H^1_P(\G_{0,i},M)$ and $(h_2)_i\in 
H_P^1(\G_{0,i}^{\fp},M)$ similarly. Equality (\ref{sofar}) 
translates to 
\be \label{sofargc} -g_i|_{\G_{1,i}^{\fp}} = (h_1)_i|_{\G_{1,i}^{\fp}} = 
(h_2)_i|_{\G_{1,i}^{\fp}} \ee
Fix $g_{\textup{mod},i} \in \phi^{-1}_{\G_{1,i}^{\fp}}(g)$ and regard it 
as an element of $\Hom_{\mathbf{Z}}(D_0,M)$ invariant under $\G_{0,i}$. 
Using diagram 
(\ref{diagraminclusion0}) with $\G_i = \G_{1,i}^{\fp}$ and $\G'_i 
= \G_{0,i}^{\fp}$ and equality (\ref{sofargc}) we conclude that there 
exists $(h_2)_{\textup{mod},i} \in 
\phi^{-1}_{\G_{0,i}^{\fp}}\bigl((h_2)_i\bigr)$ such that 
$(h_2)_{\textup{mod},i}  = -g_{\textup{mod},i}$ regarded as elements of 
$\Hom_{\mathbf{Z}}(D_0,M)$. Hence $g_{\textup{mod},i}$ is invariant under 
both $\G_{0,i}$ and $\G_{0,i}^{\fp}$. 

\begin{lemma} \label{lemmageneration} For $i=1, \dots , \# \Cl_F$ the
groups $\G_{0,i}$ and $\G_{0,i}^{\fp}$ generate $\G_{-1,i}$. \end{lemma}

\begin{proof} This is an immediate consequence of Theorem 3 on page 110 in
\cite{Serre77}.  \end{proof}

Using Lemma \ref{lemmageneration} 
we conclude that 
$g_{\textup{mod},i} \in \Hom_{\mathbf{Z}[\G_{-1,i}]}(D_0,M)$. Put $g'_i = 
\phi_{\G_{-1,i}}(g_{\textup{mod},i})$. Again, by the commutativity of 
diagram (\ref{diagraminclusion0}) with $\G'_i = \G_{-1,i}$ and $\G_i = 
\G_{0,i}$ we have $g'_i|_{\G_{0,i}} = g_i$. Hence $g':=(g'_i)_i \in 
\bigoplus_i H^1_P(\G_{-1,i},M) \cong H^1_!(X_{-1}, \tilde{M}_{-1})$ 
satisfies $g=\a_1^{-1,0}g'$. 

Thus $0=\a_1^{0,1}f + \a_{\fp}^1 \a_1^{0,1\fp} g = \a_1^{0,1} f + 
\a_{\fp}^1 \a_1^{0,1\fp} \a_1^{-1,0} g'$. By the vertical
commutativity of diagram (\ref{diagramcaug}) we have $\a_{\fp}^1 
\a_1^{0,1\fp} \a_1^{-1,0} = \a_1^{0,1} \a_{\fp}^0 \a_1^{-1,0\fp}$, so 
$0=\a_1^{0,1}f + \a_1^{0,1} \a_{\fp}^0 \a_1^{-1,0\fp} g'$. Since 
$\a_1^{0,1}$ is injective by Lemma \ref{lemmainjective}, this implies that 
$f=-\a_{\fp}^0 \a_1^{-1,0\fp} 
g'$. Hence $(f,g) = (-\a_{\fp}^0 \a_1^{-1,0\fp} g', \a_1^{-1,0} g') \in 
\beta(H^1_!(X_{-1}, \tilde{M}_{-1}))$, completing the proof of 
Proposition \ref{theorem2}.

\section{Proof of Proposition \ref{theorem3}}\label{Proof2}

In this section we prove that for a principal ideal $\fl=(l) \subset \OF$ 
such that $l \equiv 1$ (mod $\fN$) 
we have $T_{\fl} f=(N\fl+1)f$ on elements $f\in H^1_!(X_{-1}, 
\tilde{M}_{-1}) \cong \bigoplus_i 
H^1_P(\G_{-1,i},M)$. For such an ideal $\fl$, the operators 
$T_{\fl}$ preserve each direct summand $H^1_P(\G_{-1,i},M)$. The 
restriction of $T_{\fl}$ to $H^1_P(\G_{-1,i},M)$ is 
given by the usual action of the double coset $\G_{-1,i} \bsmat 1\\ &l 
\esmat \G_{-1,i}$ on group cohomology (see, e.g., \cite{Hida93}). For 
$k,l\in \OF$ we put
$\sigma_{k,l}:=\bsmat 1&k\\ &l\esmat$, $\sigma_l := \bsmat l\\ &1 \esmat$. 
To describe the action of $T_{\fl}$ explicitly we use the following 
lemma.

\begin{lemma} \label{lemmadecomposition} Let $\fl=(l)$ be a principal
ideal of $\OF$ and $n \geq -1$. Then
$$\G_{n,i} \bmat1\\&l\emat \G_{n,i} = \coprod_{\substack{k\in
R(\OF/\fl)\\k\in \fp_i}}
\G_{n,i} \sigma_{k,l}
\sqcup \G_{n,i} \sigma_l,$$
where $R(\OF/\fl)$ denotes a set of representatives of $\OF/\fl$ in
$\OF$. 
\end{lemma}

\begin{proof} This is easy.  \end{proof}

\begin{lemma} \label{lemmachebotarev} Let $n\geq 3$ be an odd integer. 
Every ideal class $c$ of $F$ contains infinitely many prime ideals $\fq$ 
such 
that 
$(N\fq-1,n)=1$.

\end{lemma}

\begin{proof} We assume $F \cap \quf(\zeta_n) = \quf$, the other case 
being easier. Let $G=\Gal(F/\quf)$, $N=\Gal(\quf(\zeta_n)/\quf)$ and 
$C=\Gal(H/F)\cong \Cl_F$, where $H$ denotes the Hilbert class field of 
$F$. We have the following diagram of fields
\be\xymatrix@!C{&&F_n H \ar@{-}[dl]^C \ar@{-}[dr]^N\\
&F_n=F\quf(\zeta_n)\ar@{-}[dl]^G \ar@{-}[dr]^N & & H\ar@{-}[dl]^C\\
\quf(\zeta_n)\ar@{-}[dr]^N& & F\ar@{-}[dl]^G\\
&\quf}\ee
Choose $(\sigma,\tau) \in \Gal(F_n H/H) \times \Gal(F_n H/F_n) \cong N 
\times C$, such that $$\sigma\in \Gal(F_n H/H)\cong 
\left(\mathbf{Z}/n\mathbf{Z}\right)^{\times}$$ corresponds to an element 
$\tilde{\sigma} \in \left(\mathbf{Z}/n\mathbf{Z}\right)^{\times}$ with 
$\tilde{\sigma} \not\equiv 1$ modulo any of the divisors of $n$, and $\tau 
\in \Gal(F_n H/F_n) \cong C \cong \Cl_F$ corresponds to the ideal class 
$c$. By the Chebotarev density theorem there exist infinitely many primes 
$\mathfrak{Q}$ of the ring of integers of $F_n H$ such that 
$\Frob_{\mathfrak{Q}}=(\sigma,\tau)$. Then the infinite set of primes 
$\fq$ of 
$\OF$ lying under such $\mathfrak{Q}$ satisfy the condition of the 
lemma, i.e., $\fq \in c$ and $(N\fq-1,n)=1$.  \end{proof}

By Lemma \ref{lemmachebotarev} we may assume that the ideals 
$\fp_i$ were chosen so that 
$N\fp_i -1$ is relatively prime to the exponent of $M$ for all $i=1, 
\dots, \# \Cl_F$.

\begin{proof}[Proof of Proposition \ref{theorem3}] Let $f \in 
H^1_P(\G_{-1,i},M)$ and let $\fl=(l)$ be a principal ideal of 
$\Op$ with $l\equiv 1$ (mod $\fN$). We will prove that $f_{\fl}:= T_{\fl} 
f -(N\fl +1)=0$. By the definition of parabolic cohomology, we have 
$f_{\fl}(\G_{-1,i}\cap B)=0$ for all $B \in \mathcal{B}$. Moreover, as 
the exponent of $M$ is relatively prime to $\#\OF^{\times}$, it is 
enough to prove that $f_{\fl}(\tilde{\G}_i)=0$, where $\tilde{\G}_i:= 
\G_{-1,i} \cap \SL_2(F)$. Put 
$$\SL_2(\Op)_i:=\left\{\bmat a&b\\c&d \emat \in \SL_2(F) \mid a,d \in \Op, 
b\in \fp_i \Op, c\in \fp_i^{-1}\Op\right\}.$$ 
We first show that $f_{\fl}=0$ on the $i$-th principal congruence 
subgroup $$\G_{\fN, i}:= \left\{\bmat a&b\\c&d \emat \in \SL_2(\Op)_i \mid
b,c\in \fN \Op, a\equiv d \equiv 1 \hf \textup{mod} \hs \fN 
\Op
\right\}.$$
If $x\in \fN\fp_i\Op$, and $g \in
\SL_2(\Op)$, then $g \bmat 1&x\\&1 \emat
g^{-1} \in \SL_2(\Op) \cap \G_{-1,i}$ and $f_{\fl}\left(g \bmat 
1&x\\&1
\emat
g^{-1}\right)=0$ by the definition of parabolic cohomology. 
So $f_{\fl}=0$ on the smallest 
normal subgroup $H$ of 
$\SL_2(\Op)$ containing matrices of the form $\bmat 1&x\\ &1 \emat$ with 
$x\in \fN \fp_i \Op$. By a 
theorem of Serre \cite{Serre70}, $$H=\G_{\fN\fp_i} := \left\{\bmat 
a&b\\c&d \emat \in 
\SL_2(\Op) \mid 
b,c\in \fN\fp_i \Op, a\equiv d \equiv 1 \hf \textup{mod} \hs \fN \fp_i \Op 
\right\}.$$ Thus $f_{\fl}=0$ on $\G_{\fN\fp_i}$. Put
$$\G'_{\fN\fp_i}:= \left\{\bmat a&b\\c&d \emat \in \SL_2(\Op)_i \mid
b,c\in \fN \fp_i \Op
\right\}.$$
Since $\G_{\fN\fp_i}$ is a normal subgroup of $\G'_{\fN\fp_i}$ of 
index $N\fp_i -1$ we have $f_{\fl} \left(\G'_{\fN\fp_i}\right)=0$ by Lemma 
\ref{lemmachebotarev} and our choice of $\fp_i$. On one 
hand $f_{\fl}$ is zero on the elements of the form $\bmat 1\\ c &1 \emat$, 
$c \in \fN \fp_i^{-1} \Op$, (again by the definition of parabolic 
cohomology) and on the other hand elements of this form together with 
$\G'_{\fN\fp_i}$ generate $\G_{\fN,i}$, so $f_{\fl}(\G_{\fN,i})=0$, as 
asserted. 

Thus $f_{\fl}$ descends to the quotient $\tilde{\G}_i/\G_{\fN,i}$. 
However, 
on this quotient all 
$\sigma_{k,l}$ and $\sigma_l$ act as the identity, since $l\equiv 
1$ 
(mod $\fN$) and we can always choose $k\in \fN\fp_i$. Thus $f_{\fl}=0$.  
 \end{proof}

\bibliography{standard1}

\end{document}